\documentclass[12pt]{article}

\newcommand{\nd}{\noindent}

\newcommand{\vu}{\vspace{.1cm}}

\newcommand{\vd}{\vspace{.2cm}}

\newcommand{\vt}{\vspace{.3cm}}

\newcommand{\vq}{\vspace{.4cm}}

\newcommand{\qed}{\nolinebreak\hfill$\Box$\par\medbreak}

\newcommand{\F}{\mathcal{F}}

\usepackage{amsmath,amsthm,amssymb}

\usepackage{amsxtra,amscd, amsbsy}

\usepackage{eucal}

\begin{document}

\begin{center}
\Large{\bf Partial groupoid actions: globalization, Morita theory
and Galois theory}
\end{center}

\begin{center}
{\bf Dirceu Bagio$^1$ and Antonio Paques$^2$}\\

\vt

{ \footnotesize
$^{1}$Departamento de Matem\'atica\\
Universidade Federal de Santa Maria\\
97105-900, Santa Maria, RS, Brazil \\
E-mail: {{\it bagio@smail.ufsm.br} }

\vu

$^{2}$Instituto de Matem\' atica\\
Universidade Federal do Rio Grande do Sul\\
91509-900, Porto Alegre, RS, Brazil\\
E-mail: {\it paques@mat.ufrgs.br}}
\end{center}

\vt

\begin{abstract}
\nd{\scriptsize In this paper we introduce the notion of a partial
action of a groupoid on a ring as well as we give a criteria for the
existence of a globalization of it. We construct a Morita context
associated to a globalizable partial groupoid action and we
introduce the notion of a partial Galois extension, which is related
to the strictness of this context.}
\end{abstract}



\vt

\vt

\vq

\nd{\large{\bf Introduction}}

\vd

The notion of a partial groupoid action is introduced in this paper
as a generalization of a partial group action \cite{E98, DE} as well
as of a partial ordered groupoid action \cite{BFP}. The notion of a
partial group action was introduced by Exel in \cite{E98} and it was
motivated by questions originated in the context of group actions on
$C^*$-algebras (see, for instance, \cite{AEE}, \cite{E94} and
\cite{McC}). A purely algebraic study of partial group actions was
done in \cite{DE} and a corresponding Galois theory was developed in
\cite{DFP}, stimulating further investigations by several others
authors, see \cite{D} for a more extensive bibliography.

\vu

Partial group actions can be easily obtained as restrictions of
global ones, and this fact led to the important question on knowing
under what conditions a given partial action is of this type. This
question was first considered by Abadie \cite{Ab} in the context of
continuous partial group actions on topological spaces and
$C^\star$-algebras. The algebraic version of globalization (or
enveloping action) of a partial group action was given by Dokuchaev
and Exel in \cite{DE}. A nice approach on the relevance of the
relationship between partial and global group actions, in several
branches of mathematics, can be seen in \cite{D}.

\vu

Groups are particular examples of groupoids (or ordered groupoids)
and any partial group action has, as in the global case, a natural
groupoid associated (see \cite{A}, \cite{Gil} and \cite{B}). Thus,
it is natural to consider the study and the development of an
algebraic theory on partial groupoid actions.

\vu

Partial ordered groupoid actions were introduced by Gilbert in
\cite{Gil} as ordered premorphisms and a study of them from an
algebraic point of view was presented in \cite{BFP}.

\vu

We start this paper by introducing the formal definition of a
partial action $\alpha$ of a groupoid $G$ on a ring $R$ (section 1)
and the formal definition of a globalization of $\alpha$ as a global
action $\beta$ of $G$ on a suitable ring $T$, satisfying some
appropriate conditions (section 2). We also give necessary and
sufficient conditions for the existence of such a globalization
(section 2), generalizing the corresponding result of \cite{DE}, as
well as we prove that the partial skew groupoid ring $R\star_\alpha
G$ (whose notion was introduced in \cite{BFP}) and the (global) skew
groupoid ring $T\star_\beta G$ are Morita equivalent, whenever the
existence of $(T,\beta)$ is ensured (section 3). This previous result
and its proof are the versions for partial groupoid actions of the
corresponding adaptations to a purely algebraic setting made in
\cite[Theorem 5.2]{DE} of \cite[Theorem 3.3]{Ab}.

\vu

For the rest of the paper is assumed that $\alpha$ is globalizable,
with globalization $(T,\beta)$. In the section 4 we introduce the
notion of the subring $R^\alpha$ of the elements of $R$ that are
invariant under the partial action $\alpha$, we prove that the rings
$R^\alpha$ and $T^\beta$ are isomorphic, and we construct a Morita
context relating the rings $R^\alpha$ and $R\star_\alpha G$. In the
section 5 we present the notion of an $\alpha$-partial Galois
extension of $R$, associating to it the Morita context constructed
in the section 4.

An action of a groupoid $G$ on a $K$-algebra $R$ ($K$ being a
subring of the center of $R$) naturally induces an structure of a
$KG$-module algebra on $R$, in the sense of \cite{CG}. In the
section 6 we briefly discuss the relationship between these two
notions.

\vu

Throughout this paper by ring we mean a non-necessarily commutative
and non-necessarily unital ring.

\vq

\vt

\nd{\large{\bf 1. Partial actions of groupoids}}

\vd

In this section we will introduce the notion of a partial groupoid
action, which is slightly more general than the notion of a partial
ordered groupoid action given in \cite{BFP}. We start recalling some
notions and notations related to groupoids.

\vu

Groupoids are usually presented as small categories in which every
morphism is invertible. But they can also be regarded, as algebraic
structures, as a natural generalization of groups. We will adopt
here the algebraic version of a groupoid given in \cite{L}. A {\it
groupoid} is a non-empty set $G$ equipped with a partially defined
binary operation, that we will denote by concatenation, for which
the usual axioms of a group hold whenever they make sense, that is:

\vu

(i) For every $g,h,l\in G$, $g(hl)$ exists if and only if $(gh)l$
exists and in this case they are equal.

\vu

(ii) For every $g,h,l\in G$, $g(hl)$ exists if and only if $gh$ and
$hl$ exist.

\vu

(iii) For each $g\in G$ there exist (unique) elements $d(g),r(g)\in
G$ such that $gd(g)$ and $r(g)g$ exist and $gd(g)=g=r(g)g$.

\vu

(iv) For each $g\in G$ there exists an element $g^{-1}\in G$ such
that $d(g)=g^{-1}g$ and $r(g)=gg^{-1}$.

\vu

It follows from this definition that the element $g^{-1}$ is unique
with the properties described in (iv) as well as that
$(g^{-1})^{-1}=g$, for every $g\in G$. Furthermore, for every
$g,h\in G$, the element $gh$ exists if and only if $d(g)=r(h)$ if
and only if $h^{-1}g^{-1}$ exists and, in this case, we have
$(gh)^{-1}=h^{-1}g^{-1}$, $d(gh)=d(h)$ and $r(gh)=r(g)$. We will
denote by $G^2$ the subset of the pairs $(g,h)\in G\times G$ such
that $d(g)=r(h)$.

\vu

An element $e\in G$ is called an {\it identity} of $G$ if
$e=d(g)=r(g^{-1})$, for some $g\in G$. In this case $e$ is called
the {\it domain identity} of $g$ and the {\it range identity} of
$g^{-1}$. We will denote by $G_0$ the set of all identities of $G$.
For any $e\in G_0$ we have  $d(e)=r(e)=e$, $e^{-1}=e$ and we will
denote by $G_e$ the set of all $g\in G$ such that $d(g)=r(g)=e$.
Clearly, $G_e$ is a group called the {\it isotropy (or principal)
group associated to $e$}.

\vu

A {\it partial action $\alpha$} of a groupoid $G$ on a ring $R$ is a
pair  $$\alpha=(\{D_g\}_{g\in G},\{\alpha_g\}_{g\in G}),$$ where for
each $g\in G$, $D_{r(g)}$ is an ideal of $R$, $D_g$ is an ideal of
$D_{r(g)}$ and $\alpha_g:D_{g^{-1}}\longrightarrow D_g$ is an
isomorphism of rings, and the following conditions hold:

\vu

(i) $\,\,\alpha_e$ is the identity map $I_{D_e}$ of $D_e$,

\vu

(ii) $\,\,\alpha_h^{-1}(D_{g^{-1}}\cap D_h)\subset D_{(gh)^{-1}}$,

\vu

(iii) $\,\,\alpha_g(\alpha_h(x))=\alpha_{gh}(x)$, for every $x\in
\alpha_h^{-1}(D_{g^{-1}}\cap D_h)$,

\vu

\nd for all $e\in G_0$ and $(g,h)\in G^2$.

\vu

The domain (resp., range) of the composition $\alpha_g\alpha_h$ is,
by definition, the largest domain (resp., range) where it makes
sense, that is, the domain of $\alpha_g\alpha_h$ is given by
$\text{dom}(\alpha_g\alpha_h)=\alpha_h^{-1}(D_{g^{-1}}\cap D_h)$ and
its range by $\text{ran}(\alpha_g\alpha_h)=\alpha_g(D_{g^{-1}}\cap
D_h)$.

\vu

Note that, by the conditions (ii) and (iii), the map $\alpha_{gh}$
is an extension of the composition $\alpha_g\alpha_h$. We say that
$\alpha$ is {\it global} if $\alpha_g\alpha_h=\alpha_{gh}$ for all
$(g,h)\in G^2$.

\vt

\nd{\bf Lemma 1.1}\,\,{\sl  Let $\alpha=\left(\{D_g\}_{g\in
G},\{\alpha_g\}_{g\in G}\right)$ be a partial action of a groupoid
$G$ on a ring $R$. Then, the following statements hold:}

\vu

(i) {\sl $\alpha$ is global if only if $D_g=D_{r(g)}$ for all $g\in
G$.}

\vu

(ii) {\sl $\alpha_{g}^{-1}=\alpha_{g^{-1}}$, for all $g\in G$.}

\vu

(iii) {\sl $\alpha_g(D_{g^{-1}}\cap D_h)=D_g\cap D_{gh}$, for all
$(g,h)\in G^2$.}

\vd

\nd{\bf Proof.}\,\, Suppose that  $\alpha$ is a global action. So,
$\alpha_g(D_{g^{-1}}\cap
D_h)=\text{ran}(\alpha_g\alpha_h)=\text{ran}(\alpha_{gh})=D_{gh},$
for all $(g,h)\in G^2$. In particular, $D_g=D_{r(g)}$, for all $g\in
G$. Conversely, assume that $D_g=D_{r(g)}$ for all $g\in G$. Then,

\[
\begin{array}{ccl}
\text{dom}(\alpha_g\alpha_h)&=&\alpha_h^{-1}(D_{g^{-1}}\cap D_h)=\alpha_h^{-1}(D_{d(g)}\cap D_{r(h)})\\
&=&\alpha_h^{-1}(D_{r(h)})=\alpha_h^{-1}(D_h)=D_{h^{-1}}\\
&=&D_{d(h)}=D_{d(gh)}=D_{(gh)^{-1}}\\
&=&\text{dom}(\alpha_{gh}).
\end{array}
\]
for all $(g,h)\in G^2$. Consequently,
$\alpha_g\alpha_h=\alpha_{gh}$.

\vu

The assertion (ii) is an immediate consequence from the definition
of a partial action. The last statement has a proof similar to that
of \cite[Corollary 2.2 (ii)]{BFP}. \qed

\vd

\nd{\bf Remark 1.2}\,\, It is immediate to see from the above that
the pair $\alpha_{(e)}=(\{D_g\}_{g\in G_e}, \{\alpha_g\}_{g\in
G_e})$ is a partial action (in the sense of \cite{DE}) of the group
$G_e$ on the ring $D_e$, for every $e\in G_0$.

\vt

\nd{\bf Example 1.3}\,\, We can obtain examples of partial action of
a groupoid by restriction of a global action, in a standard way.
Indeed, consider a global action $\beta=(\{E_g\}_{g\in G},
\{\beta_g\}_{g\in G})$ of a groupoid $G$ on a ring $T$ and for each
$e\in G_0$ let $D_e$ be an ideal of $E_e$ (for instance, if $I$ is any
ideal of $T$ it suffices to take $D_e=I\cap E_e$). Then,
$D_g=D_{r(g)}\cap\beta_g(D_{d(g)})$ is an ideal of $D_{r(g)}$ and
$\alpha_g=\beta_g|_{D_{g^{-1}}}$ is an isomorphism of rings, for all
$g\in G$. Now, take $R=\prod_{e\in G_0}D_e$ and, for each $e\in
G_0$, let $\iota_e:D_e\to R$ denote the map given by
$\iota_e(x)=(x_l)_{l\in G_0}$, with $x_e=x$ and $x_l=0$ for all
$l\neq e$. Setting $\iota_g=\iota_{r(g)}|_{D_g}$,
$D'_g=\iota_g(D_g)$ and
$\alpha'_g=\iota_g\alpha_g\iota_{g^{-1}}^{-1}:D'_{g^{-1}}\to D'_g$,
it is immediate to verify that $\alpha'=(\{D'_g\}_{g\in
G},\{\alpha'_g\}_{g\in G})$ is a partial action of $G$ on $R$. If,
in particular, each $D_e$, $e\in G_0$, is of the type $D_e=I\cap
E_e$ for some ideal $I$ of $T$, then $\alpha=(\{D_g\}_{g\in
G},\{\alpha_g\}_{g\in G})$ is also a partial action of $G$ on $I$
(see \cite[Example 2.6(2)]{BFP}). We shall say that $\alpha'$
(resp., $\alpha$) is a {\it restriction} of $\beta$ to $R$ (resp.,
$I$).

\vd

Now, take (in the above example) the ideal
$\overline{E}_g=\sum_{{r(h)=r(g)}}\beta_h(D_{d(h)})$ of $E_g$ and
set $\overline{\beta}_g=\beta_g|_{\overline{E}_{g^{-1}}}$, for each
$g\in G$. Let $\overline{T}=\prod_{e\in G_0}\overline{E}_e$ and
$\iota_e:\overline{E}_e\to\overline{T}$ be the map defined similarly
as above, for every $e\in G_0$. Setting,
$E'_g=\iota_{r(g)}(\overline{E}_g)$ and
$\beta'_g=\iota_{r(g)}\overline{\beta}_g\iota_{d(g)}^{-1}:E'_{g^{-1}}\to
E'_g$, one can easily see that $R\subseteq\overline{T}$,
$D'_g\subseteq E'_g$ and $\beta'=(\{E'_g\}_{g\in G},
\{\beta'_g\}_{g\in G})$ is also a global action of $G$ on
$\overline{T}$, whose restriction to $R$ also coincides with
$\alpha'$. We will see, in the next section, that any partial action
$\alpha=(\{D_g\}_{g\in G},\{\alpha_g\}_{g\in G})$ of a groupoid $G$
on a ring $R$ can be obtained in this way ``up to some equivalence",
provided that $D_g$ is a unital ring for all $g\in G$.

\vt

\vq

\nd{\large{\bf 2. Globalization}}

\vd

A global action $\beta=\left(\{E_g\}_{g\in G},\{\beta_g\}_{g\in
G}\right)$ of a groupoid $G$ on a ring $T$ is called a {\it
globalization} of a partial action $\alpha=\left(\{D_g\}_{g\in
G},\{\alpha_g\}_{g\in G}\right)$ of $G$ on a ring $R$ if, for each
$e\in G_0$, there is a monomorphism  of rings
$\varphi_{e}:D_{e}\rightarrow E_{e}$ such that the following
properties are satisfied:

\vu

(i) $\varphi_{e}(D_e)$ is an ideal of $E_{e}$,

\vu

(ii) $\varphi_{r(g)}(D_g)=\varphi_{r(g)}(D_{r(g)})\cap
\beta_g(\varphi_{d(g)}(D_{d(g)}))$,

\vu

(iii) $\beta_g\circ\varphi_{d(g)}(a)=\varphi_{r(g)}\circ
\alpha_g(a)$, for all $a\in D_{g^{-1}}$,

\vu

(iv) $E_g=\sum_{r(h)=r(g)}\beta_h(\varphi_{d(h)}(D_{d(h)}))$.

\vu

We say that $\beta$ is {\it unique up to equivalence} if for any
global action  $\beta'=\left(\{E'_g\}_{g\in G},\{\beta'_g\}_{g\in
G}\right)$ of $G$ on a ring $T'$, which also is a globalization of
$\alpha$, there exists, for each $e\in G_0$, an isomorphism of rings
$\psi_{e}:E'_{e}\to E_{e}$ such that
$\beta_g\circ\psi_{d(g)}(a)=\psi_{r(g)}\circ\beta'_g(a)$, for all
$a\in E'_{d(g)}$.

\vt

\nd{\bf Theorem 2.1}\,\,{\sl  Let $\alpha=\left(\{D_g\}_{g\in
G},\{\alpha_g\}_{g\in G}\right)$ be a partial action of a groupoid
$G$ on a ring $R$ and suppose that $D_e$ is a unital ring for each
$e\in G_0$. Then, $\alpha$ admits a globalization $\beta$ if and
only if each ideal $D_g$, $g\in G$, is a unital ring. Furthermore,
if $\beta$ exists then it is unique up to equivalence.}

\vd

\nd{\bf Proof.}\,\, If $\beta=\left(\{E_g\}_{g\in
G},\{\beta_g\}_{g\in G}\right)$ is a globalization for $\alpha$ and
$\varphi_{e}:D_{e}\rightarrow E_{e}$, $e\in G_0$, are the
corresponding ring monomorphisms, then
$\varphi_{r(g)}(D_g)=\varphi_{r(g)}(D_{r(g)})\cap\beta_g(\varphi_{d(g)}(D_{d(g)}))$
is clearly a unital ring so $D_g$ also is, for every $g\in G$.

\vu

Conversely, assume that each $D_g$, $g\in G$, is a unital ring with
identity $1_g$. Thus, $1_g$ is a central idempotent of $R$ and
$D_g=D_{r(g)}1_g=R1_g$.

Let $\F=\F(G,R)$ be the ring of all maps from $G$ into $R$ and, for
each $g\in G$, put $X_g=\{h\in G\,\,:\,\,r(h)=r(g)\}$ and
$F_g=\{f\in\F\ |\ f(h)=0,\text{for all}\ h\not\in X_g\}$. Clearly,
$F_g$ is an ideal of $\F$ and $F_g=F_{r(g)}$. From now on and
according to the notational convenience, we will also denote the
value $f(h)\in R$ by $f|_h$, for all $f\in F$ and $h\in G$.

For $g\in G$ and $f\in F_{g^{-1}}$ let
$\beta_g:F_{g^{-1}}\rightarrow F_g$ be the map given by

$$\beta_g(f)|_h=\left\{\begin{array}{l}
f(g^{-1}h),\,\,\,\text{if}\,\,\, h\in X_g  \\
0, \,\,\,\text{otherwise}\,\,\,
\end{array}\right.$$

Notice that if $h\in X_g$ then $r(h)=r(g)=d(g^{-1})$ which implies
that the product $g^{-1}h$ exists and $r(g^{-1}h)=r(g^{-1})$. Thus,
$\beta_g$ is well defined. Clearly, $\beta_g$ is a ring
homomorphism. Furthermore, for $g\in G$, $f\in F_{g^{-1}}$ and $h\in
X_{g^{-1}}$ we have

\[
\begin{array}{ccl}
\beta_{g^{-1}}\circ\beta_g(f)|_h &=& \beta_g(f)|_{gh}=f(g^{-1}(gh))\\
&=& f(d(g)h)=f(r(h)h)=f(h).\
\end{array}
\]
Similarly, $\beta_g\circ\beta_{g^{-1}}(f)|_h=f(h)$, for all $f\in
F_g$ and $h\in X_g$. Thus, $\beta_g$ is a ring isomorphism.

Note that $\beta_e=I_{F_e}$, for all $e\in G_0$ and
$\beta_{gh}(f)|_l=f((h^{-1}g^{-1})l)=f(h^{-1}(g^{-1}l))=\beta_g\circ\beta_h(f)|_l$,
for all $f\in F_{h^{-1}}$, $l\in X_g$ and $(g,h)\in G^2$. Hence,
$\beta=(\{F_g\}_{g\in G}, \{\beta_g\}_{g\in G})$ is a global action
of $G$ on $\F$.

\vu

Now, for each $e\in G_0$, define $\varphi_e:D_e\rightarrow F_e$
given by
$$\varphi_e(a)|_h=\left\{\begin{array}{l}
\alpha_{h^{-1}}(a1_h),\,\,\,\text{if}\,\,\, r(h)=e  \\
0, \,\,\,\text{otherwise}
\end{array}\right.$$
for all $a\in D_e$ and $h\in G$. It follows directly from the
definition that $\varphi_e(a)|_e=a$. Thus, $\varphi_e$ is a
monomorphism of rings, for all $e\in G_0$.

\vu

Let $E_g$ be the subring of $F_g$ generated by
$\bigcup_{r(h)=r(g)}\beta_h(\varphi_{d(h)}(D_{d(h)}))$, for all
$g\in G$. Notice that $\varphi_{d(g)}(D_{d(g)})\subseteq E_{d(g)}$.
Let $T=\prod_{e\in G_0}E_e$ and, for each $e\in G_0$, let
$\iota_e:E_e\to T$ be the injective map given by
$\iota_e(x)=(x_l)_{l\in G_0}$, with $x_e=x$ and $x_l=0$ for all
$l\neq e$. For convenience of notation we will identify $E_e$ with
$\iota_e(E_e)$ and $\varphi_e$ with $\iota_e\circ\varphi_e$, as well
as we will denote also by the same $\beta_g$, the ring isomorphism
given by the composition map
$\iota_{r(g)}\circ\beta_g|_{E_{g^{-1}}}\circ\iota_{d(g)}^{-1}$ from
$\iota_{d(g)}(E_{g^{-1}})\equiv E_{g^{-1}}$ onto
$E_{g}\equiv\iota_{r(g)}(E_g)$. By construction,
$\beta=(\{E_g\}_{g\in G},\{\beta_g\}_{g\in G})$ is a global action
of $G$ on $T$. Our goal is to show that $\beta$ is a globalization
of $\alpha$. We start by checking the property (iii) of the
definition of a globalization, namely:
$$\beta_g\circ\varphi_{d(g)}(a)=\varphi_{r(g)}\circ
\alpha_g(a),\,\,\,\text{ for all}\,\,\, a\in D_{g^{-1}}.$$

Consider $g\in G$, $a\in D_{g^{-1}}$ and $h\in G$. There are two
possibilities to consider. In the first, $r(h)=r(g)$. In this case,
$\varphi_{r(g)}(\alpha_g(a))|_h=\alpha_{h^{-1}}(\alpha_g(a)1_h)$ and
$\beta_g(\varphi_{d(g)}(a))|_h=\varphi_{d(g)}(a)|_{g^{-1}h}=\alpha_{h^{-1}g}(a1_{g^{-1}h})$.
However, by \cite[Corollary 2.2]{BFP},
$\alpha_{h^{-1}g}(a1_{g^{-1}h})\in \alpha_{h^{-1}g}(D_{g^{-1}}\cap
D_{g^{-1}h})\subset D_{h^{-1}}$ and
$\alpha_{h^{-1}g}(1_{g^{-1}}1_{g^{-1}h})=1_{h^{-1}}1_{h^{-1}g}$.
Since $a1_{g^{-1}h}1_{g^{-1}}\in D_{g^{-1}}\cap D_{g^{-1}h}$, it
follows that
\[
\begin{array}{ccl}
\beta_g(\varphi_{d(g)}(a))|_h &=& \alpha_{h^{-1}g}(a1_{g^{-1}h})=\alpha_{h^{-1}g}(a1_{g^{-1}h})1_{h^{-1}}\\
&=&\alpha_{h^{-1}g}(a1_{g^{-1}h})\alpha_{h^{-1}g}(1_{g^{-1}}1_{g^{-1}h})=
\alpha_{h^{-1}g}(a1_{g^{-1}h}1_{g^{-1}})\\
&=&\alpha_{h^{-1}}(\alpha_g(a1_{g^{-1}h}1_{g^{-1}}))=
\alpha_{h^{-1}}(\alpha_g(a1_{g^{-1}})\alpha_g(1_{g^{-1}}1_{g^{-1}h}))\\
&=&\alpha_{h^{-1}}(\alpha_g(a)1_h)= \varphi_{r(g)}(\alpha_g(a))|_h.\
\end{array}
\]
In the second, $r(h)\neq r(g)$. In this case, we have
$\varphi_{r(g)}(\alpha_g(a))|_h=0=\beta_g(\varphi_{d(g)}(a))|_h$.

\vd

The next step is to show that
$$
\varphi_{r(g)}(D_g)=\varphi_{r(g)}(D_{r(g)})\cap
\beta_g(\varphi_{d(g)}(D_{d(g)})).
$$

Given $g\in G$ and $c\in \varphi_{r(g)}(D_{r(g)})\cap
\beta_g(\varphi_{d(g)}(D_{d(g)}))$, there exist $a\in D_{r(g)}$ and
$b\in D_{d(g)}$ such that
$c=\varphi_{r(g)}(a)=\beta_g(\varphi_{d(g)}(b))$. Then,
$a=a1_{r(g)}=\alpha_{r(g)^{-1}}(a1_{r(g)})=\varphi_{r(g)}(a)|_{r(g)}=\beta_g(\varphi_{d(g)}(b))|_{r(g)}=
\varphi_{d(g)}(b)|_{g^{-1}r(g)}=\varphi_{d(g)}(b)|_{g^{-1}}=
\alpha_g(b1_{g^{-1}})\in D_g$ and so $c=\varphi_{r(g)}(a)\in
\varphi_{r(g)}(D_g)$. Conversely, if $c\in \varphi_{r(g)}(D_g)$ then
$c=\varphi_{r(g)}(a)$, for some $a\in D_g$. Taking
$b=\alpha_{g^{-1}}(a)\in D_{g^{-1}}$ we have
$\beta_g(\varphi_{d(g)}(b))\overset{\text{(iii)}}{=}\varphi_{r(g)}(\alpha_g(b))=
\varphi_{r(g)}(a)$. So, $c\in\varphi_{r(g)}(D_{r(g)})\cap
\beta_g(\varphi_{d(g)}(D_{d(g)}))$.

\vd

To conclude that $\beta$ is a globalization of $\alpha$, it remains
to prove that $\varphi_e(D_e)$ is an ideal of $E_e$ for all $e\in
G_0$. To see this it is enough to check that
$\varphi_{e}(b)\beta_h(\varphi_{d(h)}(a))$ and
$\beta_h(\varphi_{d(h)}(a))\varphi_{e}(b)$ are elements of
$\varphi_{e}(D_{e})$, for all $h\in X_e$, $a\in D_{d(h)}$ and $b\in
D_{e}$. Given $l\in G$, we have again two possibilities to consider.
In the first, $r(l)=r(h)=e$. Since $\alpha_{h}(a1_{h^{-1}})b\in
D_{e}$, we obtain
\[
\begin{array}{ccl}
\beta_h(\varphi_{d(h)}(a))\varphi_{e}(b)|_l &=&
(\varphi_{d(h)}(a)|_{h^{-1}l})(\varphi_{e}(b)|_{l})=
\alpha_{l^{-1}h}(a1_{h^{-1}l})\alpha_{l^{-1}}(b1_l)\\
&=&\alpha_{l^{-1}}(\alpha_h(a1_{h^{-1}})b1_l)=\varphi_{e}(\alpha_h(a1_{h^{-1}})b)|_l.\
\end{array}
\]
In the second, $r(l)\neq r(h)=e$ and so
$$\beta_h(\varphi_{d(h)}(a))\varphi_{e}(b)|_l=0=\varphi_{e}(\alpha_h(a1_{h^{-1}})b)|_l.$$
Hence,
$\beta_h(\varphi_{d(h)}(a))\varphi_{e}(b)=\varphi_{e}(\alpha_h(a1_{h^{-1}})b)\in
\varphi_{e}(D_{e})$. Similarly we also have
$\varphi_{e}(b)\beta_h(\varphi_{d(h)}(a))\in \varphi_{e}(D_{e})$.

\vu

To end the proof it is required to show the uniqueness (up to
equivalence) of the globalization $\beta$ of $\alpha.$ Notice that
$E_g=\sum_{r(h)=r(g)}\beta_h(\varphi_{d(h)}(D_{d(h)}))$, for every
$g\in G$. This is an immediate consequence from the fact that each
$\varphi_{d(h)}(D_{d(h)})$ is an ideal of $E_{d(h)}$, as proved
above.

\vu

Now, suppose that $\beta'=(\{E'_g\}_{g\in G}, \{\beta'_g\}_{\in G})$
is a global action of $G$ on a ring $T'$ and also a globalization of
$\alpha$. Then, for each $e\in G_0$, there exists a rind
monomorphism $\varphi'_{e}:D_{e}\rightarrow E'_{e}$ satisfying the
conditions (i)-(iv) of the definition of a globalization.

For each $e\in G_0$, consider the map $\eta_{e}:E'_e\rightarrow E_e$
given by
$$\sum_{1\leq i\leq n}\beta'_{h_i}(\varphi'_{d(h_i)}(a_i))\mapsto \sum_{1\leq i\leq n}\beta_{h_i}
(\varphi_{d(h_i)}(a_i)),$$ with $h_i\in X_e$ and $a_i\in
D_{d(h_i)}$, for all $1\leq i\leq n$. First of all, we need to check
that $\eta_{e}$ is well defined.

Suppose that $\sum_{1\leq i\leq
n}\beta'_{h_i}(\varphi'_{d(h_i)}(a_i))=0$. For all $l\in X_e$, with
$d(l)=e$, and $r\in D_{d(l)}$, we have $\sum_{1\leq i\leq
n}\beta'_{h_i}(\varphi'_{d(h_i)}(a_i))\beta'_{l}(\varphi'_{d(l)}(r))=0,$
which implies that $\sum_{1\leq i\leq
n}\beta'_{l^{-1}h_i}(\varphi'_{d(h_i)}(a_i))\varphi'_{d(l)}(r)=0.$
Since
$\varphi'_{d(l^{-1}h_i)}(D_{d(l^{-1}h_i)})=\varphi'_{d(h_i)}(D_{d(h_i)})$
is an ideal of $E'_{d(h_i)}$, it follows that
$\beta'_{l^{-1}h_i}(\varphi'_{d(l^{-1}h_i)}(D_{d(h_i)}))$ is an
ideal of $E'_{r(l^{-1})}$. Note that
$\beta'_{l^{-1}h_i}(\varphi'_{d(h_i)}(a_i))\varphi'_{d(l)}(r)$ is
contained in the following ideal
\begin{eqnarray*}
&&\beta'_{l^{-1}h_i}(\varphi'_{d(l^{-1}h_i)}(D_{d(h_i)}))\cap
\varphi'_{d(l)}(D_{d(l)})=\\
&&\beta'_{l^{-1}h_i}(\varphi'_{d(l^{-1}h_i)}(D_{d(l^{-1}h_i)}))\cap
\varphi'_{r(l^{-1}h_i)}(D_{r(l^{-1}h_i)})\overset{\text{(ii)}}{=}\\
&&\varphi'_{r(l^{-1}h_i)}(D_{l^{-1}h_i})=
\varphi'_{d(l)}(D_{d(l)}1_{l^{-1}h_i})=\varphi'_{d(l)}(D_{d(l)})\varphi'_{d(l)}(1_{l^{-1}h_i})
\end{eqnarray*}
whose identity element is $\varphi'_{d(l)}(1_{l^{-1}h_i})$.
Therefore,
\[
\begin{array}{ccl}
\beta'_{l^{-1}h_i}(\varphi'_{d(h_i)}(a_i))\varphi'_{d(l)}(r)&=&
\beta'_{l^{-1}h_i}(\varphi'_{d(h_i)}(a_i))\varphi'_{d(l)}(1_{l^{-1}h_i})\varphi'_{d(l)}(r)\\
&=&\beta'_{l^{-1}h_i}(\varphi'_{d(h_i)}(a_i))\varphi'_{d(l)}
(\alpha_{l^{-1}h_i}(1_{h_i^{-1}l}))\varphi'_{d(l)}(r)\\
&\overset{\text{(iii)}}{=}&\beta'_{l^{-1}h_i}(\varphi'_{d(h_i)}(a_i))\beta'_{l^{-1}h_i}
(\varphi'_{d(h_i)}(1_{h_i^{-1}l}))\varphi'_{d(l)}(r)\\
&=&\beta'_{l^{-1}h_i}\circ\varphi'_{d(h_i)}(a_i1_{h_i^{-1}l})\varphi'_{d(l)}(r)\\
&\overset{\text{(iii)}}{=}&\varphi'_{d(l)}(\alpha_{l^{-1}h_i}(a_i1_{h_i^{-1}l})r).
\end{array}
\]
In a similar way we also get
\begin{align*}
\beta_{l^{-1}h_i}(\varphi_{d(h_i)}(a_i))\varphi_{d(l)}(r)=
\varphi_{d(l)}(\alpha_{l^{-1}h_i}(a_i1_{h_i^{-1}l})r).
\end{align*}
Thus, $$0=\sum_{1\leq i\leq
n}\beta'_{l^{-1}h_i}(\varphi'_{d(h_i)}(a_i))\varphi'_{d(l)}(r)=
\sum_{1\leq i\leq
n}\varphi'_{d(l)}(\alpha_{l^{-1}h_i}(a_i1_{h_i^{-1}l})r).$$ Since
$\varphi'_{d(l)}$ is a monomorphism, it follows that $\sum_{1\leq
i\leq n}\alpha_{l^{-1}h_i}(a_i1_{h_i^{-1}l})r=0$. Hence,
$$0=\sum_{1\leq i\leq n}\varphi_{d(l)}(\alpha_{l^{-1}h_i}(a_i1_{h_i^{-1}l})r)=
\sum_{1\leq i\leq
n}\beta_{l^{-1}h_i}(\varphi_{d(h_i)}(a_i))\varphi_{d(l)}(r)$$ and
applying $\beta_l$ we obtain
$$\sum_{1\leq i\leq n}\beta_{h_i}(\varphi_{d(h_i)}(a_i))\beta_l(\varphi_{d(l)}(r))=0, \quad\text{for all}\quad
r\in D_{d(l)}.$$ Hence, the element $a=\sum_{1\leq i\leq n}
\beta_{h_i}(\varphi_{d(h_i)}(a_i))$ annihilates
$\beta_l(\varphi_{d(l)}(D_{d(l)}))$, for all $l\in G$ with $r(l)=e$.
In particular, $a$ annihilates the ideal $L$ of $E_g$ given by
$\sum_{1\leq i\leq n}\beta_{h_i}(\varphi_{d(h_i)}(D_{d(h_i)}))$,
which is unital by \cite[Lemma 4.4]{DE}. Since $a\in L$, it follows
that $a=0$ and so $\eta_{e}$ is a well defined homomorphism of
rings. Actually, $\eta_{e}$ is an isomorphism, whose inverse is
given by $\sum_{1\leq i\leq
n}\beta_{h_i}(\varphi_{d(h_i)}(a_i))\mapsto \sum_{1\leq i\leq
n}\beta'_{h_i}(\varphi'_{d(h_i)}(a_i)),$ with $r(h_i)=e$ and $a_i\in
D_{d(h_i)}$, for all $1\leq i\leq n$.

\vu

Finally, for each $u=\beta'_h(\varphi'_{d(h)}(a))$, with
$r(h)=r(d(g))=d(g)$, and $a\in D_{d(h)}$ we have
$\eta_{r(g)}\circ\beta'_g(u)=\eta_{r(g)}(\beta'_{gh}(\varphi'_{d(h)}(a)))=
\beta_{gh}(\varphi_{d(h)})(a)=\beta_g(\beta_h(\varphi_{d(h)}(a)))=\beta_g\circ\eta_{d(g)}(u)$,
which completes the proof.\qed

\vt

\nd{\bf Example 2.2}\,\,\,Let $G=\{d(g),r(g),g,g^{-1}\}$ be a
groupoid and $R=Ke_1\oplus Ke_2\oplus Ke_3$, where $K$ is a unital
ring and $e_1, e_2, e_3$ are pairwise orthogonal central idempotents
with sum $1_R$. Take $D_{d(g)}=Ke_1\oplus Ke_2$,
$D_{r(g)}=D_g=Ke_3$, $D_{g^{-1}}=Ke_1$,
$\alpha_{d(g)}=I_{D_{d(g)}}$, $\alpha_{r(g)}=I_{D_{r(g)}}$,
$\alpha_g(ae_1)=ae_3$, $\alpha_{g^{-1}}(ae_3)=ae_1$ and note that
$\alpha=(\{D_g\}_{g\in G},\{\alpha_g\}_{g\in G})$ is a partial
action of $G$ on $R$.

\vspace{0.1cm}

Observe that $X_{g^{-1}}=\{d(g),g^{-1}\}$, $X_{g}=\{r(g),g\}$,
$F_{g^{-1}}\simeq R\times R$ via $f\mapsto (f(d(g)),f(g^{-1}))$ and
$F_g\simeq R\times R$ via $f\mapsto (f(r(g)),f(g))$. Assuming the
natural identifications as in the proof of Theorem 2.1, we have that
$\beta_g:F_{g^{-1}}\rightarrow F_g,\,(r,s)\mapsto(s,r)$, for all
$(r,s)\in R\times R$. Also, $\varphi_{d(g)}(D_{d(g)})\simeq
\{(a,\alpha_g(a1_{g^{-1}}))\,:\,a\in D_{d(g)}\}$ and
$\varphi_{r(g)}(D_{r(g)})\simeq
\{(b,\alpha_{g^{-1}}(b1_{g}))\,:\,b\in D_{r(g)}\}$. Since $b1_g=b$,
for all $b\in D_{r(g)}$, it follows that
$\beta_{g^{-1}}(\varphi_{r(g)}(D_{r(g)}))\subset
\varphi_{d(g)}(D_{d(g)})$. Consequently,
$E_{d(g)}=E_{g^{-1}}=\varphi_{d(g)}(D_{d(g)})\simeq D_{d(g)}$.

\vspace{0.1cm}

On the other hand, we have that
$\beta_{g}(\varphi_{d(g)}(D_{d(g)}))(e_3,e_1)\subset
\varphi_{r(g)}(D_{r(g)})$. Then,
\[
\begin{array}{ccl}
E_{r(g)}&=&E_g\\
&=&\varphi_{r(g)}(D_{r(g)})\oplus \beta_{g}(\varphi_{d(g)}(D_{d(g)}))(1_R-e_3,1_R-e_1)\\
&=& D_{r(g)}\oplus K
\end{array}
\]
Hence, the globalization of $\alpha$ is the action
$\beta=(\{E_g\}_{g\in G},\{\beta_g\}_{g\in G})$ of $G$ on $T=R\oplus
Ke_4=Ke_1\oplus Ke_2\oplus Ke_3\oplus Ke_4$, where $e_1, e_2, e_3,
e_4$ are pairwise orthogonal central idempotents with sum $1_T$,
$\beta_{d(g)}=I_{E_{d(g)}}$, $\beta_{r(g)}=I_{E_{r(g)}}$,
$\beta_g(ae_1+be_2)=ae_3+be_4$ and
$\beta_{g^{-1}}(ae_3+be_4)=ae_1+be_2$, for all $a,b\in K$.

\vt

\nd{\bf Remark 2.3}\,\, Let $\alpha=\left(\{D_g\}_{g\in
G},\{\alpha_g\}_{g\in G}\right)$ be a partial action of a groupoid
$G$ on a ring $R$ having a globalization $\beta=\left(\{E_g\}_{g\in
G},\{\beta_g\}_{g\in G}\right)$. Simplifying notation, assume that
$D_e\subset E_e$, for all $e\in G_0$. Note that the group $G_e$ acts
globally on $E_e=\sum_{r(h)=e}\beta_h(D_{d(h)})$, and it is clear
that the action of $G_e$ on the subalgebra $\sum_{h\in
G_e}\beta_h(D_{e})$ of $E_e$ is a globalization (in the sense of
\cite{DE}) of the partial action $\alpha_{(e)}$ of $G_e$ on $D_e$.
In particular, Theorem 2.1 is a generalization of \cite[Theorem
4.5]{DE}.

\vt

\vq

\nd{\large{\bf 3. Morita Equivalence}}

\vd

Let $\alpha=(\{D_g\}_{g\in G}, \{\alpha_g\}_{g\in G})$ be a partial
action of a groupoid $G$ on a ring $R$ and assume that each $D_g$ is
unital. Let $\beta=(\{E_g\}_{g\in G}, \{\beta_g\}_{g\in G})$ be a
global action of $G$ on a ring $T$ and a globalization of $\alpha$.
Simplifying notation, we will assume that $D_{e}\subseteq E_{e}$,
for all $e\in G_0$.

In this section we will see, under suitable conditions, that the
partial skew groupoid ring $R\star_{\alpha} G$ and the skew groupoid
ring $T\star_{\beta}G$ are Morita equivalent.

\vu

We start recalling that a {\it Morita context} is a six-tuple
$(R,R',M,M',\tau,\tau')$ where $R$ and $R'$ are unital rings, $M$ is
an $(R,R')$-bimodule, $M'$ is an $(R',R)$-bimodule and
$\tau:M\otimes_{R'}M'\rightarrow R$ and
$\tau':M'\otimes_{R}M\rightarrow R'$ are bimodule maps satisfying
the following conditions:

\vu

(i) $\tau(x\otimes x')y=x\tau'(x'\otimes y)$, for all $x,y\in M$ and
$x'\in M'$,

\vu

(ii) $\tau'(x'\otimes x)y'=x'\tau(x\otimes y')$, for all $x',y'\in
M'$ and $x\in M$.

\vu

\nd Following \cite[Theorems 4.1.4 and 4.1.17]{Rw}, if $\tau$ and
$\tau'$ are surjective then the categories of left $R$-modules and
left $R'$-modules are equivalent and the rings $R$ and $R'$ are
called {\it Morita equivalent}. One also says, in this case, that
the corresponding Morita context is {\it strict}.

\vd

Following \cite[Section 3]{BFP}, the partial skew groupoid ring
$R\star_\alpha G$ corresponding to $\alpha$ is defined as the direct
sum
$$R\star_\alpha G=\bigoplus_{g\in G}D_g\delta_g$$ in which the
$\delta_g$'s are symbols, with the usual addition, and
multiplication determined by the rule
\[
(x\delta_g)(y\delta_h)=
\begin{cases}
\alpha_g(\alpha_{g^{-1}}(x)y)\delta_{gh} &\text{if $(g,h)\in G^2$}\\
0  &\text{otherwise},
\end{cases}
\]
for all $g, h\in G$, $x\in D_g$ and $y\in D_h$.

\vu

Let $A$ (resp. $B$) denote the partial skew groupoid ring (resp.,
the skew groupoid ring ) $R\star_{\alpha}G$ (resp.,
$T\star_{\beta}G$). Let also $1_g$ denote the identity element of
$D_g$, for all $g\in G$. As already observed in the former section,
$1_g$ is a central idempotent of $R$ and it is immediate to check
that it also is a central idempotent of $T$. So, we also have
$D_g=D_{r(g)}1_g=E_{r(g)}1_g=R1_g=T1_g$, for all $g\in G$.
Furthermore, it follows from $D_g=D_{r(g)}\cap\beta_g(D_{d(g)})$
that $1_g=1_{r(g)}\beta_g(1_{d(g)})$. We will also assume that $G_0$
is finite, which is equivalent to say that $A$ is unital  with
$1_A=\sum_{e\in G_0}1_e\delta_e$, by \cite[Proposition 3.3]{BFP}.

\vt

\nd{\bf Proposition 3.1}\,\,{\sl Let $R$, $T$, $G$, $\alpha$,
$\beta$, $A$ and $B$ be as above. Then:}

\vu

(i) {\sl $B1_A=\left\{\sum_{g\in
G}c_g\delta_g\,\,:\,\,c_g\in\beta_g(D_{d(g)}),\,\,\,\text{for
all}\,\,\,g\in G\right\}$,}

\vu

(ii) {\sl $1_AB=\left\{\sum_{g\in G}c_g\delta_g\,\,:\,\,c_g\in
D_{r(g)},\,\,\,\text{for all}\,\,\,g\in G\right\}$,}

\vu

(iii) {\sl $1_AB1_A=A$,}

\vu

(iv) {\sl $B1_AB=B$. }

\vd

\nd{\bf Proof.}\,\, (i) For every $g\in G$ and $s\in E_g$ we have
$s=\beta_g(t)$, for some $t\in E_{g^{-1}}=E_{d(g)}$, and
$$(s\delta_g)1_A=(s\delta_g)(1_{d(g)}\delta_{d(g)})=
s\beta_g(1_{d(g)})\delta_g=\beta_g(t1_{d(g)})\delta_g$$ with
$\beta_g(t1_{d(g)})\in \beta_g(D_{d(g)})$, since $D_{d(g)}$ is an
ideal of $E_{d(g)}$.

For the reverse inclusion, take $a\in D_{d(g)}$ and
$c_g=\beta_g(a)$. Then,
$$c_g\delta_g=\beta_g(a)\delta_g=\beta_g(a1_{d(g)})\delta_g=
(c_g\delta_g)(1_{d(g)}\delta_{d(g)})=(c_g\delta_g)1_A.$$

\vu

(ii) For $s\in E_g$ we have
$$1_A(s\delta_g)=(1_{r(g)}\delta_{r(g)})(s\delta_g)=1_{r(g)}\beta_{r(g)}(s)\delta_g=
1_{r(g)}s\delta_g$$ with $1_{r(g)}s \in D_{r(g)}$, because
$D_{r(g)}$ is an ideal of $E_g$.

Conversely, if $c_g\in D_{r(g)}$ then
$$c_g\delta_g=(1_{r(g)}\delta_{r(g)})(c_g\delta_g)=1_A(c_g\delta_g).$$

\vu

(iii) Consider $a\in E_g$. Then, $a1_g\in D_g$ and
\[
\begin{array}{ccl}
1_A(a\delta_g)1_A&=&(1_{r(g)}a\delta_g)(1_{d(g)}\delta_{d(g)})\\
&=&1_{r(g)}a\beta_g(1_{d(g)})\delta_g\\
&=&1_ga\delta_g.\
\end{array}
\]
The converse is immediate.

\vu

(iv) It is enough to show that $B\subseteq B1_AB$. Since, for each
$h\in G$, $E_h=\sum_{r(g)=r(h)}\beta_g(D_{d(g)})$, the result
follows from (i) and (ii) because
$$\beta_g(a)\delta_h=(\beta_g(a)\delta_g)(1_{r(g^{-1}h)}\delta_{g^{-1}h})\in (B1_A)(1_AB)=B1_AB,$$
for all $a\in D_{d(g)}$. \qed

\vt

\nd{\bf Theorem 3.2}\,\,{\sl Let $R$, $T$, $G$, $\alpha$, $\beta$,
$A$ and $B$ be as above and suppose that $B$ is unital. Then
$A=R\star_{\alpha}G$ and $B=T\star_{\beta}G$ are Morita equivalent.}

\vd

\nd{\bf Proof.}\,\, Clearly, $M=1_AB$ is a $(A,B)$-bimodule and
$N=B1_A$ is a $(B,A)$-bimodule. By Proposition 3.1, we have that
$\tau:M\otimes_BN\rightarrow A,\,\,\,m\otimes n\mapsto mn$ and
$\tau':N\otimes_AM\rightarrow B,\,\,\,n\otimes m\mapsto nm$ are
surjective bimodule maps and one can easily see by a straightforward
calculation that the six-tuple $(A,B,M,N,\tau,\tau')$ is a Morita
context. \qed

\vt

\vq

\nd{\large{\bf 4. The subring of invariants and the trace map}}

\vd

Throughout the rest of this paper $\alpha=(\{D_g\}_{g\in G},
\{\alpha_g\}_{g\in G})$ will denote a partial action of a finite
groupoid $G$ on a ring $R$ such that each $D_g$ is unital with
identity element $1_g$. Recall that, in this case, $1_g$ is a
central idempotent of $R$ and $D_g=D_{r(g)}1_g=R1_g$, for all $g\in
G$.

We assume that $\beta=\left(\{E_g\}_{g\in G},\{\beta_g\}_{g\in
G}\right)$ is a globalization of $\alpha$, acting on a ring $T$.
Recalling the construction of a globalization given in section 2,
and in order to simplify notation, we also assume that
$T=\bigoplus_{e\in G_0}E_e$, $R=\bigoplus_{e\in G_0}D_e$ and $D_e$
is an ideal of $E_e$, for all $e\in G_0$. Clearly, $R$ is a unital
ring with identity element given by $1_R=\sum_{e\in G_0}1_e$ and,
under the assumptions considered, $R=T1_R.$

\vu

Two relevant concepts which appear in Galois theory are the notions
of subring of invariants and trace map. The subring of {\it
invariants} of $R$ under $\alpha$ is defined similarly as in
\cite{DFP} by
$$R^{\alpha}=\{x\in R\,:\,\alpha_g(x1_{g^{-1}})=x1_g,\,\,\,\text{for all}\,\,\, g\in G\}.$$

\nd{\bf Remark 4.1}\,\,\, Notice that $R^\alpha\subseteq
\bigoplus_{e\in G_0}D_e^{\alpha_{(e)}}.$ Indeed, any $x\in R$ is of
the form $x=\sum_{e\in G_0}x_e$, with $x_e\in D_e$, and $x\in
R^\alpha$ if and only if $\alpha_g(x_{d(g)}1_{g^{-1}})=x_{r(g)}1_g$,
for every $g\in G$. Similarly, we have that $a=\sum_{e\in G_0}a_e\in
T^{\beta}$ if only if $\beta_g(a_{d(g)})=a_{r(g)}$, for all $g\in
G$. In general, the inclusion $R^\alpha\subseteq \bigoplus_{e\in
G_0}D_e^{\alpha_{(e)}}$ is  strict as it is well shown in the
following example: $G=\{d(g),r(g),g,g^{-1}\}$, $R=Ke_1\oplus
Ke_2\oplus Ke_3\oplus Ke_4\oplus Ke_5$, where $K$ is a unital ring
and $e_1, e_2, e_3, e_4, e_5$ are pairwise orthogonal central
idempotents with sum $1_R$, $D_{d(g)}=D_{g^{-1}}=Ke_1\oplus Ke_2$,
$D_{r(g)}=Ke_3\oplus Ke_4\oplus Ke_5$, $D_g=Ke_3\oplus Ke_4$,
$\alpha_{d(g)}=I_{D_{d(g)}}$, $\alpha_{r(g)}=I_{D_{r(g)}}$,
$\alpha_g(ae_1+be_2)=ae_3+be_4$ and
$\alpha_{g^{-1}}(ce_3+de_4)=ce_1+de_2$, for all $a,b,c,d\in K$. It
is immediate to check that $\alpha=(\{D_g\}_{g\in
G},\{\alpha_g\}_{g\in G})$ is a partial (not global) action of $G$
on $R$, $R^{\alpha}=K(e_1+e_3)\oplus K(e_2+e_4)\oplus Ke_5$,
$D_{d(g)}^{\alpha_{(d(g))}}=Ke_1\oplus Ke_2$ and
$D_{r(g)}^{\alpha_{(r(g))}}=Ke_3\oplus Ke_4\oplus Ke_5$.

\vt

The trace map is defined as the map  $t_\alpha:\,R\rightarrow R$
given by
$$t_\alpha(x)=\sum_{g\in G}\alpha_g(x1_{g^{-1}}),$$ for all $x\in R$.

\vt

\nd{\bf Lemma 4.2}\,\,\,{\sl The map $t_\alpha$ is a homomorphism of
$(R^\alpha,R^\alpha)$-bimodules and}

(i) {\sl  $t_\alpha(R)\subseteq R^{\alpha}$,}

\vu

(ii) {\sl $t_\alpha(\alpha_g(x))=t_\alpha(x)$, for all $x\in
D_{g^{-1}}$ and $g\in G$.}

\vd

\nd{\bf Proof.}\,\ The first assertion is obvious.

(i) Consider $x\in R$ and $h\in G$ and observe that
$\alpha_g(x1_{g^{-1}})1_{h^{-1}}\in D_g\cap D_{h^{-1}}$. Then,
\[
\begin{array}{ccl}
\alpha_h(t_\alpha(x)1_{h^{-1}})&=&\sum_{g\in
G}\alpha_h(\alpha_g(x1_{g^{-1}})1_{h^{-1}})\\
&=&\sum_{r(g)=d(h)}\alpha_{hg}(x1_{(hg)^{-1}})1_h\\
&=&\sum_{r(l)=r(h)}\alpha_l(x1_{l^{-1}})1_h\\
&=&\sum_{l\in G}\alpha_l(x1_{l^{-1}})1_h\\
&=&t_\alpha(x)1_h.\
\end{array}
\]

\vt

(ii) Recall from \cite[Corollary 2.2]{BFP} that, for all $g,h\in G$
such that $d(g)=r(h)$, we have $\alpha_g(D_{g^{-1}}\cap D_h)=D_g\cap
D_{gh}$ and consequently $\alpha_g(1_{g^{-1}}1_h)=1_g1_{gh}$. Thus,
\[
\begin{array}{ccl}
t_\alpha(\alpha_g(x))&=&\sum_{h\in
G}\alpha_h(\alpha_g(x)1_{h^{-1}})\\
&=&\sum_{d(h)=r(g)}\alpha_{hg}(x1_{(hg)^{-1}})1_h\\
&=&\sum_{d(l)=d(g)}\alpha_l(x1_{l^{-1}})1_{lg^{-1}}\\
&=&\sum_{d(l)=d(g)}\alpha_l(x1_{l^{-1}})\alpha_l(1_{l^{-1}}1_{g^{-1}})\\
&=&\sum_{d(l)=d(g)}\alpha_l(x1_{l^{-1}})=t_\alpha(x)\
\end{array}
\]
for all $x\in D_{g^{-1}}$ and $g\in G$. \qed

\vt

\nd{\bf Remark 4.3}\,\, The statement (i) in Lemma 4.2 is of
fundamental importance to introduce a notion of a Galois extension
for partial groupoid actions, which generalizes that defined for
partial group actions in \cite{DFP}. And this statement holds
because the assumptions assumed on the ring $R$ in the beginning of
this section. In  general, Lemma 4.2(i)-(ii) is not true, as we show
in the following example.

Consider the groupoid $G=\{g_1,g_2,g_3\}$ with $G_0=\{g_1,g_2\}$,
$g_3^{-1}=g_3$ and $g_3g_3=g_2$, and take $R=Ke_1\oplus Ke_2\oplus
Ke_3\oplus Ke_4$, where $K$ is a unital ring and $e_1, e_2, e_3,
e_4$ are pairwise orthogonal central idempotents with sum $1_R$. Put
$D_{g_1}=R$, $D_{g_2}=Ke_2\oplus Ke_3\oplus Ke_4$,
$D_{g_3}=Ke_2\oplus Ke_3$, $\alpha_{g_1}=I_{D_{g_1}}$,
$\alpha_{g_2}=I_{D_{g_2}}$ and $\alpha_{g_3}(ae_2+be_3)=be_2+ae_3$,
for all $a,b\in K$. By a simple calculation, we have that
$\alpha=(\{D_g\}_{g\in G},\{\alpha_g\}_{g\in G})$ is a partial (not
global) action of $G$ on $R$ and it is immediate to check that
$R^\alpha=Ke_1\oplus K(e_2+e_3)\oplus Ke_4$,
$t_\alpha(e_3)=e_2+2e_3\not\in R^\alpha$ and $t_\alpha(e_1)=e_1\neq
0=t_\alpha(\alpha_{g_3}(e_11_{g_3^{-1}}))$.

\vt

Now we observe that by a straightforward and routine calculation it
follows that $R$ is a $(R^{\alpha}, R\star_{\alpha} G)$-bimodule
(resp., $(R\star_{\alpha} G,R^{\alpha})$-bimodule), with the right
(resp., left) action of $R\star_{\alpha} G$ on $R$ given by $x\cdot
a\delta_g=\alpha_{g^{-1}}(xa)$ (resp., $a\delta_g\cdot x=
a\alpha_g(x1_{g^{-1}})$), for all $g\in G$, $a\in D_g$ and $x\in R$.

\vu

Furthermore, the map $R\times R\rightarrow R^{\alpha},\,\,
(x,y)\mapsto t_\alpha(xy)$ (resp., $R\times R\rightarrow
R\star_{\alpha} G,\,\, (x,y)\mapsto\sum_{g\in
G}x\alpha_g(y1_{g^{-1}})\delta_g$) is $R\star_{\alpha} G$-balanced
(resp., $R^{\alpha}$-balanced) by Lemma 4.2. So we can consider the
maps
$$\tau\,:\,R\otimes_{R\star_{\alpha} G}R\rightarrow R^{\alpha},\,\,\tau(x\otimes y) =
t_\alpha(xy)$$ and $$\tau'\,:\,R\otimes_{R^{\alpha}}R\rightarrow
R\star_{\alpha} G,\,\,\tau'(x\otimes y)=\sum_{g\in
G}x\alpha_g(y1_{g^{-1}})\delta_g.$$

\vd

\nd{\bf Proposition 4.4}\,\,{\sl
$(R\star_{\alpha}G,R^{\alpha},R,R,\tau, \tau')$ is a Morita context.
The map $\tau$ is surjective if and only if there exists $x\in R$
such that $t_\alpha(x)=1_R$.}

\vd

\nd{\bf Proof.}\,\  The first assertion follows by a straightforward
and routine calculation and the second is obvious.\qed

\vt

Our goal in the rest of this section is to prove that
$R^{\alpha}\simeq T^{\beta}$. It is clear that $T^{\beta}1_R\subset
R^{\alpha}$. In fact, for $a=\sum_{e\in G_0}a_e\in T^{\beta}$ and
$g\in G$, we have

\[
\begin{array}{ccl}
\alpha_g((a1_R)1_{g^{-1}})&=&\alpha_g(a_{d(g)}1_{d(g)}1_{g^{-1}})=\beta_g(a_{d(g)}1_{g^{-1}})\\
&=&\beta_g(a_{d(g)})\beta_g(1_{g^{-1}})=a_{r(g)}1_g=(a1_R)1_g.\
\end{array}
\]

For each $g\in G$, we recall that $X_g=\{h\in G\,\,:\,\,r(h)=r(g)\}$
(see section 2). Now, assume that
$X_e=\{g_{1,e}=e,g_{2,e},\ldots,g_{n_e,e}\}$, for every $e\in G_0$.
Then, $E_e=\sum_{1\leq j\leq n_e}\beta_{g_{j,e}}(D_{d(g_{j,e})})$ is
a unital ring \cite[Lemma 4.4]{DE} and its identity element $1'_e$
is given by the following boolean sum $$1'_e=\sum_{1\leq l\leq
n_e}\sum_{i_1\leq\ldots\leq
i_l}(-1)^{l+1}\beta_{g_{i_1,e}}(1_{d(g_{i_1,e})})\ldots\beta_{g_{i_l,e}}(1_{d(g_{i_l,e})}),$$
of the central idempotents $\beta_{g_{j,e}}(1_{d(g_{j,e})})$ of
$E_e$, $1\leq j\leq n_e$. We can also write $1'_e$ as the orthogonal
sum $1'_e=v_{1,e}\oplus v_{2,e}\oplus\ldots\oplus v_{n_e,e},$ where
$v_{1,e}=1_e$ and
$$v_{j,e}=(1'_e-1_e)(1'_e-\beta_{g_{2,e}}(1_{d(g_{2,e})}))\ldots
(1'_e-\beta_{g_{j-1,e}}(1_{d(g_{j-1,e})}))\beta_{g_{j,e}}(1_{d(g_{j,e})}),$$
for all $2\leq j\leq n_e$. It also follows from the above that $T$
is a unital ring, with identity element $1_T=\sum_{e\in G_0}1'_e$.
Moreover, the map $\psi_e:T\longrightarrow E_e$ defined by
$$\psi_e(a)=\sum_{1\leq l\leq n_e}\sum_{i_1\leq\ldots\leq i_l}(-1)^{l+1}\beta_{g_{i_1,e}}
(1_{d(g_{i_1,e})})\ldots\beta_{g_{i_l,e}}(1_{d(g_{i_l,e})})
\beta_{g_{i_l,e}}(a_{d(g_{i_l,e})}),$$ for all $e\in G_0$ and
$a=\sum_{e\in G_0}a_e\in T$, is clearly a homomorphism of rings. It
is immediate that $\psi_e(1_R)=1'_e$ and we can also write
$\psi_e(a)=\sum_{1\leq j\leq
n_e}\beta_{g_{j,e}}(a_{d(g_{j,e})})v_{j,e}$.

\vt

\nd{\bf Lemma 4.5}\,\,{\sl If $x\in R^{\alpha}$ then
$\beta_g(\psi_{d(g)}(x))=\psi_{r(g)}(x)$, for all $g\in G$.}

\vd

\nd{\bf Proof.}\,\ Take $g,l,h\in G$ with $r(l)=r(h)$ and $x\in
R^{\alpha}$. Then,

\[
\begin{array}{ccl}
\beta_l(x1_{d(l)})\beta_h(1_{d(h)})&=&\beta_l(x1_{d(l)}\beta_{l^{-1}h}(1_{d(h)}))\\
&=&\beta_l(x1_{r(l^{-1}h)}\beta_{l^{-1}h}(1_{d(l^{-1}h)}))\\
&=&\beta_{l}(x1_{l^{-1}h})=\beta_{l}(\alpha_{l^{-1}h}(x1_{h^{-1}l}))\\
&=&\beta_l(\beta_{l^{-1}h}(x1_{r(h^{-1}l)}\beta_{h^{-1}l}(1_{d(h^{-1}l)})))\\
&=&\beta_h(x1_{d(h)}\beta_{h^{-1}l}(1_{d(l)}))\\
&=&\beta_h(x1_{d(h)})\beta_l(1_{d(l)}).\
\end{array}
\]
Since $gX_{d(g)}=X_{r(g)}$, it follows from the equality above that
$\beta_g(V)=W$, where $V$ (resp., $W$) denote the set of all
summands of $\psi_{d(g)}(x)$ (resp., $\psi_{r(g)}(x)$).
Consequently, $\beta_g(\psi_{d(g)}(x))=\psi_{r(g)}(x)$.\qed

\vd

The next result is an extension of \cite[Proposition 2.3]{DFP} to
the context of partial groupoid actions.

\vt

\nd{\bf Theorem 4.6}\,\,{\sl The subrings of invariants $R^{\alpha}$
and $T^{\beta}$ are isomorphic.}

\vd

\nd{\bf Proof.}\,\ Consider the map $\psi:T\longrightarrow T$
defined by $\psi(a)=\sum_{e\in G_0}\psi_e(a),$ for any $a\in T$.
Note that, by Lemma 4.5,
$\beta_g(\psi(x)_{d(g)})=\beta_g(\psi_{d(g)}(x))=\psi_{r(g)}(x)=\psi(x)_{r(g)}$,
for all $x\in R^{\alpha}$ and $g\in G$. Thus,
$\psi(R^{\alpha})\subset T^{\beta}$. Since $\psi_e(x)1_e=x_e1_e$ it
follows that $\psi(x)1_R=x$, for all $x=\sum_{e\in G_0}x_e\in
R^{\alpha}$. Furthermore, for $a=\sum_{e\in G_0}a_e\in T^{\beta}$,
we have

\[
\begin{array}{ccl}
\psi_e(a1_R)&=&\sum\limits_{j=1}^{n_e}\beta_{g_{j,e}}(a_{d(g_{j,e})}1_{d(g_{j,e})})v_{j,e}\\
&=&\sum\limits_{j=1}^{n_e}\beta_{g_{j,e}}(a_{d(g_{j,e})})\beta_{g_{j,e}}(1_{d(g_{j,e})})v_{j,e}\\
&=&\sum\limits_{j=1}^{n_e}\beta_{g_{j,e}}(a_{d(g_{j,e})})v_{j,e}=\sum\limits_{j=1}^{n_e}a_{r(g_{j,e})}v_{j,e}\\
&=&a_e\left(\sum\limits_{j=1}^{n_e}v_{j,e}\right)=a_e1'_e=a_e.\
\end{array}
\]
So, $\psi|_{R^{\alpha}}:R^{\alpha}\longrightarrow T^{\beta}$ is a
ring isomorphism whose inverse is the map given by the
multiplication by $1_R$.\qed

\vt

\vq

\nd{\large{\bf 5. Galois theory}}

\vd

The notion of a partial Galois extension for partial group actions
given in \cite{DFP} is a generalization of the classical notion of
Galois extension presented in \cite{CHR}. In the sequel, we
introduce this notion for partial groupoid actions. Recall that $R$,
$T$, $G$, $\alpha$ and $\beta$ are assumed here just like as in the
former section. So all the notations introduced there will be also
freely used in this section.

\vu

We will say that $R$ is an {\it $\alpha$-partial Galois extension}
of the subring of invariants $R^\alpha$ if there exist elements
$x_i,y_i\in R$, $1\leq i\leq m$, such that $\sum_{1\leq i\leq
m}x_i\alpha_g(y_i1_{g^{-1}})=\delta_{e,g}1_e$, for all $e\in G_0$
and $g\in G$. In particular, in this case, $D_e$ is an
$\alpha_{(e)}$-partial Galois extension of $D_e^{\alpha_{(e)}}$ in
the sense of \cite{DFP}, for every $e\in G_0$. The set
$\{x_i,y_i\}_{1\leq i\leq n}$ is called a {\it partial Galois
coordinate system} of $R$ over $R^\alpha$. When $\alpha$ is global,
we say simply that $R$ is an {\it $\alpha$-Galois extension of
$R^{\alpha}$} and that {\it $\{x_i,y_i\}_{1\leq i\leq n}$ is a
Galois coordinate system}.

\vt

\nd{\bf Theorem 5.1}\,\, \, {\sl The following statements hold:}

\vu

(i) {\sl If $T$ is a $\beta$-Galois extension of $T^{\beta}$ then
$R$ is an $\alpha$-partial Galois extension of $R^\alpha$ .}

\vu

(ii) {\sl Suppose that $G$ satisfies the condition:
$g^{-1}_{j,r(l)}lg_{j,d(l)}\in G_0$ implies that $l\in G_0$, for all
$1\leq j\leq n_{d(l)}=n_{r(l)}$. Then the converse of (i) also
holds.}

\vd

\nd{\bf Proof.} (i) Let $\{a_i,b_i\}_{1\leq i\leq n}$ be a Galois
coordinate system of $T$ over $T^{\beta}$ and suppose that
$a_i=\sum_{e\in G_0}a_{i,e}$, $b_i=\sum_{e\in G_0}b_{i,e}$, for all
$1\leq i\leq n$. Taking $x_i=\sum_{e\in G_0}a_{i,e}1_e=a_i1_R$ and
$y_i=\sum_{e\in G_0}b_{i,e}1_e=b_i1_R$, one can easily check that
$\{x_i,y_i\}_{1\leq i\leq n}$ is a partial Galois coordinate system
of $R$ over $R^{\alpha}$.

\vu

(ii) Let $\{x_i,y_i\}_{1\leq i\leq n}$ be a partial Galois
coordinate system of $R$ over $R^{\alpha}$ and assume that
$x_i=\sum_{e\in G_0}x_{i,e}$ and $y_i=\sum_{e\in G_0}y_{i,e}$, for
all $1\leq i\leq n$. Consider the elements $a_{i,j}=\sum_{e\in
G_0}a_{i,j,e}$ and $b_{i,j}=\sum_{e\in G_0}b_{i,j,e}$, where
$a_{i,j,e}=\beta_{g_{j,e}}(x_{i,d(g_{j,e})})v_{j,e}$ and
$b_{i,j,e}=\beta_{g_{j,e}}(y_{i,d(g_{j,e})})v_{j,e}$, for all $e\in
G_0$, $1\leq i\leq n$ and $1\leq j\leq n_e$.

\vu

Notice that
\[
\begin{array}{ccl}
\sum\limits_{i,j}a_{i,j}\beta_e(b_{i,j}1'_e)&=&\sum\limits_{i,j}\beta_{g_{j,e}}(x_{i,d(g_{j,e})})
\beta_{g_{j,e}}(y_{i,d(g_{j,e})})v_{j,e}\\
&=&\sum\limits_{j}\beta_{g_{j,e}}(\sum\limits_{i}x_{i,d(g_{j,e})}y_{i,d(g_{j,e})})v_{j,e}\\
&=&\sum\limits_{j}\beta_{g_{j,e}}(1_{d(g_{j,e})})v_{j,e}=\sum\limits_{j}v_{j,e}\\
&=&1'_e,\
\end{array}
\]
for all $e\in G_0$.

\vu

Consider $u_{j,e}=(1'_e-1_e)\ldots
(1'_e-\beta_{g_{j-1,e}}(1_{d(g_{j-1,e})}))$ and
$g_{jl}=g^{-1}_{j,r(l)}lg_{j,d(l)}$, for all $e\in G_0$ and $l\notin
G_0$. Thus, $g_{jl} \notin G_0$ and consequently
\[
\begin{array}{ccl}
\sum\limits_{i,j}a_{i,j}\beta_l(b_{i,j}1'_{l^{-1}})&=&\sum\limits_{i,j}a_{i,j,r(l)}\beta_{l}(b_{i,j,d(l)})\\
&=&\sum\limits_{i,j}\beta_{g_{j,r(l)}}(x_{i,d(g_{j,r(l)})})\beta_{lg_{j,d(l)}}(y_{i,d(g_{j,d(l)})})
v_{j,r(l)}\beta_l(v_{j,d(l)})\\
&=&\sum\limits_{j}\beta_{g_{j,r(l)}}(\sum\limits_{i}x_{i,r(g_{jl})}\beta_{g_{jl}}(y_{i,d(g_{jl})})1_{d(g_{j,r(l)})})
u_{j,r(l)}\beta_l(v_{j,d(l)})\\
&=&\sum\limits_{j}\beta_{g_{j,r(l)}}(\sum\limits_{i}x_{i,r(g_{jl})}\alpha_{g_{jl}}(y_{i,d(g_{jl})}1_{g^{-1}_{jl}}))
u_{j,r(l)}\beta_l(v_{j,d(l)})\\
&=&0.\
\end{array}
\]
Hence, the set $\{a_{i,j}, b_{i,j}\}$, as above constructed, is a
Galois system for $T$ over $T^\beta$. \qed

\vd

\nd{\bf Remark 5.2}\,\,\, When $G$ is a group, the hypothesis (ii)
of Theorem 5.1 is trivially satisfied and we recover \cite[Theorem
3.3]{DFP}. Certainly, there are groupoids, others than groups, which
also satisfy the condition stated in Theorem 5.1(ii). Take, for
instance, the groupoid $G=\{g_1, g_2, g_3\}$, with $G_0=\{g_1,
g_2\}$, $g_3^{-1}=g_3$ and $g_3g_3=g_2$.

\vt

From now on, we will denote by $j:R\star_\alpha
G\to\text{End}(R)_{R^\alpha}$ the natural map given by $j(\sum_{g\in
G}a_g\delta_g)(x)=\sum_{g\in G}a_g\alpha_g(x1_{g^{-1}})$, for all
$x\in R$. Clearly, $j$ is a well defined homomorphism of left
$R$-modules. Moreover, $j$ is also a ring homomorphism. Indeed, by a
straightforward calculation one easily gets $j(1_{R\star_\alpha
G})=I_R$ and {$j((a_g\delta_g)(b_h\delta_h))= j(a_g\delta_g)\circ
j(b_h\delta_h)$ for every $g,h\in G$ such that $d(g)=r(h)$. If
otherwise $d(g)\neq r(h)$ then $(a_g\delta_g)(b_h\delta_h)=0$ and
also $D_h\cap D_{g^{-1}}=0$ which implies that $(j(a_g\delta_g)\circ
j(b_h\delta_h))(x)=a_g\alpha_g(b_h\alpha_h(x1_{h^{-1}})1_{g^{-1}})=0$,
for every $x\in R$.

\vu

For any left $R\star_\alpha G$-module $M$, we put $$M^G=\{m\in M\ |\
(1_g\delta_g)m= 1_gm,\,\,\text{for all}\,\,g\in G\}$$ the set of the
invariants of $M$ under $G$. Note that $M$ is a left $R$-module via
the embedding $x\mapsto x1_{R\star_\alpha G}=\sum_{e\in
G_0}x1_e\delta_e$ from $R$ into $R\star_\alpha G$ and, indeed, $M^G$
is a left $R^\alpha$-submodule of $M$. In particular, since $R$ is a
left $R\star_\alpha G$-module via $j$, we have that $R^G$ coincides
with the subring $R^\alpha$.

\vd

The following theorem extends (and improves) \cite[Theorem
3.1]{CHR}, \cite[Theorem 4.1]{DFP} and \cite[Theorem 3.1]{GL} to the
setting of partial groupoid actions.

\vt

\nd{\bf Theorem 5.3}\,\, \, {\sl The following statements are
equivalent:}

\vu

(i) {\sl $R$ is an $\alpha$-partial Galois extension of $R^\alpha$.}

\vu

(ii) {\sl $R$ is a finitely generated projective right
$R^\alpha$-module and $j$ is an isomorphism of rings and left
$R$-modules.}

\vu

(iii) {\sl For every left $R\star_\alpha G$-module $M$ the map
$\mu:R\otimes_{R^\alpha} M^G\to M$, given by $\mu(x\otimes m)=xm$,
is an isomorphism of left $R$-modules.}

\vu

(iv) {\sl The map $\rho:R\otimes_{R^\alpha} R\to \prod_{g\in G}
D_g$, $x\otimes y\mapsto (x\alpha_g(y1_{g^{-1}}))_{g\in G}$, is an
isomorphism of left $R$-modules.}

\vu

(v) {\sl $RtR=R\star_\alpha G$, where $t=\sum_{g\in G}1_g\delta_g$.}

\vu

(vi) {\sl The map $\tau'$ is surjective.}

\vu

(vii) {\sl $R$ is a generator for the category of the left
$R\star_\alpha G$-modules.}

\vd

{\sl Moreover, under the assumption that at least one of the above
statements holds, then the following additional statements also are
equivalent:}

\vu

(viii) {\sl $t_\alpha(R)=R^\alpha$.}

\vu

(ix) {\sl $R$ is a generator for the category of the right
$R^\alpha$-modules.}

\vu

(x) {\sl The Morita context $(R\star_{\alpha}G,R^{\alpha},R,R,\tau,
\tau')$ is strict.}

\vd

\nd{\bf Proof.} (i)$\Rightarrow$(ii)\, Take $x_i,y_i\in R$, $1\leq
i\leq n$, such that $\sum_{1\leq i\leq
n}x_i\alpha_g(y_i1_{g^{-1}})=\delta_{e,g}1_e$, for all $e\in G_0$
and $g\in G$. A dual basis for $R$ as a projective right
$R^\alpha$-module is given by the set $\{x_i,f_i\}_{i=1}^n$, with
$f_i=t_\alpha(y_i\_)\in\text{Hom}_{R^\alpha}(R,R^\alpha)$.

\vu

For any $f\in\text{End}(R)_{R^\alpha}$, taking $w=\sum_{g\in
G}\sum_{1\leq i\leq n}f(x_i)\alpha_g(y_i1_{g^{-1}})\delta_g$, an
easy calculation gives $j(w)(x)=f(x)$, for every $x\in R$. Finally,
if $v=\sum_{g\in G}a_g\delta_g$ is such that $j(v)=0$ then
\[
\begin{array}{ccl}
0&=&\sum_{h\in G}\sum_{1\leq i\leq n}j(v)(x_i)\alpha_h(y_i1_{h^{-1}})\delta_h\\
&=&\sum_{h\in G}\sum_{1\leq i\leq n}\sum_{g\in
G}a_g\alpha_g(x_i1_{g^{-1}})
\alpha_h(y_i1_{h^{-1}})\delta_h\\
&=&\sum_{r(g)=r(h)}a_g\alpha_g(\sum_{1\leq i\leq n}x_i
\alpha_{g^{-1}h}(y_i1_{h^{-1}g})1_{g^{-1}})\delta_h\\
&=&\sum_{e\in
G_0}\sum_{r(g)=r(h)}a_g\alpha_g(\delta_{e,g^{-1}h}1_e1_{g^{-1}})\delta_h
\end{array}
\]
Now, observing that $g^{-1}h=e$ implies
$e=r(e)=r(g^{-1}h)=r(g^{-1})=d(g)$, one has
$g=gd(g)=ge=g(g^{-1}h)=(gg^{-1})h=r(g)h=r(h)h=h$ and so
\[
\begin{array}{ccl}
0&=&\sum_{e\in G_0}\sum_{r(g)=r(h)}a_g\alpha_g(\delta_{e,g^{-1}h}1_e1_{g^{-1}})\delta_h\\
&=&\sum_{g\in G}a_g\alpha_g(1_{d(g)}1_{g^{-1}})\delta_g\\
&=&\sum_{g\in G}a_g\delta_g=v.\
\end{array}
\]

\vu

(ii)$\Rightarrow$(iii)\, Let $\{x_i,f_i\}_{1\leq i\leq l}$ be a dual
basis of $R$ as a right $R^\alpha$-module, and $v_i\in R\star_\alpha
G$ such that $j(v_i)=f_i$, for all $1\leq i\leq l$. It follows from
the injectivity of $j$ that $\sum_{1\leq i\leq
m}x_iv_i=1_{R\star_\alpha G}$ as well as $(1_g\delta_g)v_i=1_gv_i$,
for all $g\in G$. Thus, the map $\nu:M\to R\otimes_{R^\alpha}M^G$
given by $\nu(m)=\sum_{1\leq i\leq l}x_i\otimes v_im$ is a
well-defined homomorphism of left $R$-modules and an easy
computation shows that $\nu$ is the inverse of $\mu$.

\vd

(iii)$\Rightarrow$(iv)\, Set $\mathfrak{F}=\{f:G\to R\, |\, f(g)\in
D_g,\,\, \text{for all}\,\, g\in G\}$. Clearly, $\mathfrak{F}$ is an
$(R,R)$-bimodule naturally isomorphic to $\prod_{g\in G}D_g$ and a
left $R\star_\alpha G$-module via the well-defined action given by
$$((a_g\delta_g)f)(h)=\left\{\begin{array}{l}
a_g\alpha_g(f(g^{-1}h)1_g^{-1}),\,\,\,\text{if}\,\,\, r(h)=r(g)  \\
0, \,\,\,\text{otherwise}
\end{array}\right.$$
for all $g,h\in G$ and $f\in\mathfrak{F}$.

Also, the map $R\to\mathfrak{F}^G$, $r\mapsto f_r$, such that
$f_r(h)=\alpha_h(r1_{h^{-1}})$, $h\in G$, is a well-defined
isomorphism of left $R^\alpha$-modules. Indeed, such a map is
clearly a homomorphism of left $R^\alpha$-modules from $R$ to
$\mathfrak{F}$. Given $g,h\in G$ and $r\in R$ we have
$((1_g\delta_g)f_r)(h)=\alpha_g(f_r(g^{-1}h)1_g^{-1})=
\alpha_g(\alpha_{g^{-1}h}(r1_{h^{-1}g})1_{g^{-1}})=\alpha_h(r1_{h^{-1}})1_g=f_r(h)1_g$,
if $r(h)=r(g)$. If $r(h)\neq r(g)$, then $D_g\cap D_h=0$ and
consequently $f_r(h)1_g=0=((1_g\delta_g)f_r)(h)$. So,
$f_r\in\mathfrak{F}^G$. If for some $r\in R$ and every $h\in G$ we
have $f_r(h)=0$, then  $r=r1_R=\sum_{e\in G_0}r1_e=\sum_{e\in
G_0}\alpha_e(r1_e)=\sum_{e\in G_0}f_r(e)=0$. Finally, given
$f\in\mathfrak{F}^G$ and taking $r=\sum_{e\in G_0}f(e)$ we have
$f_r(h)=\sum_{e\in
G_0}\alpha_h(f(e)1_{h^{-1}})=\alpha_h(f(d(h))1_{h^{-1}})=\alpha_h(f(h^{-1}h)1_{h^{-1}})
=((1_h\delta_h)f)(h)=1_hf(h)=f(h),$ for every $h\in G$.

Therefore, the map given by the composition $R\otimes_{R^\alpha}
R\simeq
R\otimes_{R^\alpha}\mathfrak{F}^G\overset{\mu}{\to}\mathfrak{F}\simeq\prod_{g\in
G}D_g$ is the claimed isomorphism $\rho$.

\vd

(iv)$\Rightarrow$(i) It follows immediately from the surjectivity of
$\rho$.

\vd

(i)$\Leftrightarrow$(v) Notice that $RtR$ is an ideal of
$R\star_\alpha G$. This easily follows from the fact that
\[
\begin{array}{ccl}
(a\delta_h)t &=& \sum_{g\in G}(a\delta_h)(1_g\delta_g)=\sum_{r(g)=d(h)}a\alpha_h(1_g1_{h^{-1}})\delta_{hg}\\
&=&
\sum_{r(g)=d(h)}a1_h1_{hg}\delta_{hg}=\sum_{r(l)=r(h)}a1_l\delta_l=at\
\end{array}
\]
and
\[
\begin{array}{ccl}
t(a\delta_h)&=&\sum_{g\in G}(1_g\delta_g)(a\delta_h)=\sum_{d(g)=r(h)}\alpha_g(a1_{g^{-1}})\delta_{gh}\\
&=&\sum_{d(l)=d(h)}\alpha_{lh^{-1}}(a1_{hl^{-1}})\delta_l=\sum_{d(l)=
d(h)}\alpha_l(\alpha_{h^{-1}}(a1_h)1_{l^{-1}})\delta_l\\
&=&t\alpha_{h^{-1}}(a1_h),\
\end{array}
\]
for any $h\in G$ and $a\in D_h$. So, $RtR=R\star_\alpha G$ if and
only if there exist elements $x_i,y_i\in R$, $1\leq i\leq n$, such
that $$\sum_{e\in G_0}1_e\delta_e=1_{R\star_\alpha G}= \sum_{1\leq
i\leq n}x_ity_i=\sum_{g\in G}(\sum_{1\leq i\leq
n}x_i\alpha_g(y_i1_{g^{-1}}))\delta_g$$ if and only if
$\{x_i,y_i\}_{i=1}^n$ is a partial Galois coordinate system of $R$
over $R^\alpha$.

\vd

(v)$\Leftrightarrow$(vi) It is enough to observe that
$\tau'(R\otimes_{R^\alpha}R)=RtR$.

\vd

(ii)$\Leftrightarrow$(vii) Let $(R^\alpha)^{op}$ denote the opposite
ring of $R^\alpha$. Now, observe that $R$ has a natural left
$(R^\alpha)^{op}$-module structure via the right multiplication,
which is compatible with its left $R\star_\alpha G$-module structure
defined via $j$, that is, $R$ is a $((R^\alpha)^{op}, R\star_\alpha
G)$-bimodule. Furthermore, the map $\theta:(R^\alpha)^{op}\to
\text{End}_{R\star_\alpha G}(R)$ given by $\theta(r)(x)=xr$, for all
$r\in (R^\alpha)^{op}$ and $x\in R$, is a well defined homomorphism
of rings having inverse given by $f\mapsto f(1_R)$, for all
$f\in\text{End}_{R\star_\alpha G}(R)$. Finally, since
$\text{End}(R)_{R^\alpha}= \text{End}_{(R^\alpha)^{op}}(R)$, the
result follows by \cite[Theorem 18.8]{W}.

\vd

(viii)$\Leftrightarrow$(ix) Assume that $t_\alpha(R)=R^\alpha$.
Thus, $t_\alpha$ is a surjective right $R^\alpha$-linear map from
$R$ to $R^\alpha$ and consequently $R_{R^\alpha}$ is a generator.

Conversely, the trace ideal
$\mathfrak{T}(R)_{R^\alpha}:=\sum_{f\in\text{Hom}(R,R^\alpha)_{R^\alpha}}f(R)$
of $R^\alpha$ equals $R^\alpha$ by \cite[Theorem 13.7]{W}. Since $j$
is an isomorphism by assumption, for every
$f\in\text{Hom}(R,R^\alpha)_{R^\alpha}\subseteq\text{End}(R)_{R^\alpha}$
there exists a unique $v=\sum_{g\in G}a_g\delta_g\in R\star_\alpha
G$ such that $j(v)=f$. Therefore, $j(v)(x)\in R^\alpha$, for any
$x\in R$ and consequently
\[
\begin{array}{ccl}
j(\sum_{r(g)=r(h)}a_g1_h\delta_g)(x)&=&(\sum_{g\in G}a_g\alpha_g(x1_{g^{-1}}))1_h\\
&=&\sum_{g\in G}\alpha_h(a_g\alpha_g(x1_{g^{-1}})1_{h^{-1}})\\
&=&\sum_{g\in G}\alpha_h(a_g1_{h^{-1}})\alpha_h(\alpha_g(x1_{g^{-1}})1_{h^{-1}})\\
&=&\sum_{r(g)=d(h)}\alpha_h(a_g1_{h^{-1}})\alpha_{hg}(x1_{(hg)^{-1}}))\\
&=&\sum_{r(g)=r(h)}\alpha_h(a_{h^{-1}g}1_{h^{-1}})\alpha_g(x1_{g^{-1}}))\\
&=&j(\sum_{r(g)=r(h)}\alpha_h(a_{h^{-1}g}1_{h^{-1}})\delta_g)(x)\
\end{array}
\]
which implies $$\sum_{r(g)=r(h)}a_g1_h\delta_g=
\sum_{r(g)=r(h)}\alpha_h(a_{h^{-1}g}1_{h^{-1}})\delta_g,$$ for any
$h\in G$. Then, $a_g1_h=\alpha_h(a_{h^{-1}g}1_{h^{-1}})$ for all
$g,h\in G$ such that $r(g)=r(h)$. In particular, for $g=h$ we have
$a_g=\alpha_g(a_{d(g)}1_{g^{-1}})$, for all $g\in G$. Now, setting
$a_f=\sum_{e\in G_0}a_e$ we have
\[
\begin{array}{ccl}
t_\alpha(a_fx)&=&\sum_{g\in G}\sum_{e\in G_0}\alpha_g(a_ex1_{g^{-1}})\\
&=&\sum_{e\in G_0}\sum_{d(g)=e}\alpha_g(a_e1_{g^{-1}})\alpha_g(x1_{g^{-1}})\\
&=&\sum_{g\in G}a_g\alpha_g(x1_{g^{-1}})\\
&=& f(x),\
\end{array}
\]
for every $x\in R$. Hence,
$R^\alpha=t_\alpha(\sum_{f\in\text{Hom}(R,R^\alpha)_{R^\alpha}}a_fR)\subseteq
t_\alpha(R)$ and the result follows.

\vd

(viii)$\Leftrightarrow$(x) It follows from Proposition 4.4.\qed

\vt

\nd{\bf Corollary 5.4} \,\, {\sl Suppose that at least one of the
following assertions holds:}

\vu

(i) {\sl $R$ is a commutative ring.}

\vu

(ii) {\sl $t_\alpha(1_R)$ is invertible in $R$.}

\vu

(iii) {\sl $|G_e|1_e$ is invertible in $D_e$ and
$t_{\alpha}|_{D_e}=t_{\alpha_{(e)}}$, for all $e\in G_0$.}

\vu

{\sl Then $R$ is an $\alpha$-partial Galois extension of $R^\alpha$
if and only if the Morita context
$(R\star_{\alpha}G,R^{\alpha},R,R,\tau, \tau')$ is strict.}

\vd

\nd{\bf Proof.}\, It is enough to show that
$\mathfrak{T}(R)_{R^\alpha}=R^\alpha$.

\vu

(i) In this case the proof is the same as in \cite[Lemma 1.6]{CHR}.
Indeed, as we saw above, if $\{x_i,y_i\}_{i=1}^n$ is a partial
Galois coordinate system for $R$ over $R^\alpha$, then
$\{x_i,t_\alpha(y_i\_)\}_{i=1}^n$ is a dual basis for $R$ as a
finitely generated projective $R^\alpha$-module and so
$R\mathfrak{T}(R)_{R^\alpha}=R$. Consequently, there exists
$c\in \mathfrak{T}(R)_{R^\alpha}$ such that $1_R=t_\alpha(c)$
(see \cite[Corollary 2.2.5]{AD} for instance) and the result follows.

\vu

(ii) Put $c=t_\alpha(1_R)^{-1}$. Then,
$t_\alpha(1_R)=t_\alpha(t_\alpha(1_R)c)=t_\alpha(1_R)t_\alpha(c)$
and so $t_\alpha(c)=1_R$, which implies the required.

\vu

(iii) It follows from Remark 1.2 and \cite[Theorem 3.1 and Corollary
3.5]{BLP} that there exists $c_e\in D_e$ such that
$t_{\alpha_{(e)}}(c_e)=1_e$, for every $e\in G_0$. Now, taking
$c=\sum_{e\in G_0}c_e$ we have
$$t_\alpha(c)=\sum_{e\in G_0}t_{\alpha}(c_e)=\sum_{e\in
G_0}t_{\alpha_{(e)}}(c_e)=\sum_{e\in G_0}1_e=1_R$$ and the required
follows.\qed

\vt

\nd{\bf Remark 5.5}\, \, The condition (iii) in Corollary 5.4 is
trivially  satisfied in the case that $G$ is a group. The following
example shows that it is also satisfied even when $G$ is not a
group. Take the groupoid $G=\{g_1,g_2,g_3\}$ with $G_0=\{g_1,g_2\}$,
$g_3^{-1}=g_3$ and $g_3g_3=g_2$. Given $R=Ke_1\oplus Ke_2\oplus
Ke_3\oplus Ke_4\oplus Ke_5$, where $K$ is a unital ring and $e_1,
e_2, e_3, e_4, e_5$ are pairwise orthogonal central idempotents with
sum $1_R$, put $D_{g_1}=Ke_1\oplus Ke_2$, $D_{g_2}=Ke_3\oplus
Ke_4\oplus Ke_5$, $D_{g_3}=Ke_3\oplus Ke_4$,
$\alpha_{g_1}=I_{D_{g_1}}$, $\alpha_{g_2}=I_{D_{g_2}}$ and
$\alpha_{g_3}(ae_3+be_4)=be_3+ae_4$, for all $a,b\in K$. Note that
$\alpha=(\{D_g\}_{g\in G},\{\alpha_g\}_{g\in G})$ is a partial
action (not global) of $G$ on $R$,
$t_{\alpha}|_{D_{g_1}}=t_{\alpha_{(g_1)}}$ and
$t_{\alpha}|_{D_{g_2}}=t_{\alpha_{(g_2)}}$.

\vt

\vq

\nd{\large{\bf 6. A final remark}}

\vd

In \cite{CG} Caenepeel and De Groot developed a Galois theory for
weak Hopf algebra actions on algebras. In particular, they
considered the situation where the weak Hopf algebra is a finite
groupoid algebra (it is well known that any groupoid algebra is a
weak Hopf algebra) and a notion of groupoid action was introduced.
Actually, this previous notion and the our's one are equivalent, which was
proved by D. Fl\^ores in her PhD thesis \cite{Fl}. More specifically
she proved the following theorem.

\vd

\nd{\bf Theorem 6.1}\, \cite[Teorema 1.2.11]{Fl} {\sl Let $G$ be a
finite groupoid, $K$ a commutative ring and $R$ a $K$-algebra. Then
the following statements are equivalent:}

\vu

(i) {\sl There exists an action $\beta=(\{E_g\}_{g\in G},
\{\beta_g\}_{g\in G})$ of $G$ on $R$ such that every $E_e$, $e\in
G_0$, is unital and $R=\bigoplus_{e\in G_0}E_e.$}

\vu

(ii) {\sl $R$ has an structure of $KG$-module algebra.}

\vd

Given the action $\beta$ of $G$ on $R$ and $KG=\bigoplus_{g\in
G}Ku_g$, the $KG$-module algebra structure of $R$ is given by the
action $u_g\cdot r=\beta_g(r1_{g^{-1}})$. Conversely, given an
action $\cdot$ of $KG$ on $R$, the corresponding action $\beta$ of
$G$ on $R$ is the pair $(\{E_g\}_{g\in G}, \{\beta_g\}_{g\in G})$,
where $E_g=R1_g$, $1_g=u_g\cdot 1_R$ and
$\beta_g(r1_{g^{-1}})=(u_g\cdot r)1_g=u_g\cdot r$, for all $r\in R$ and $g\in
G$.

Finally, the assertion (iv) of Theorem 5.3 is just the definition of
Galois extension for groupoid actions as considered in \cite{CG}.

\vt

\vd

\nd{\large{\bf Acknowledgments}}

The authors would like to thank the referee for his(her) helpful
comments.



\begin{thebibliography}{9}

\bibitem{Ab} F. Abadie, {\sl Enveloping actions and Takai duality
for partial actions}, J. Funct. Analysis 197 (2003), 14-67.

\bibitem{A} F. Abadie, {\sl On partial actions and groupoids},
Proc. AMS 132 (2003), 1037-1047.

\bibitem{AEE} B. Abadie, S. Eilers and R. Exel, {\sl Morita
equivalence and crossed products by Hilbert $C^\star$-bimodules},
Trans. Amer. Math. Soc. 350 (1998), no. 8, 3043-3054.

\bibitem{AD} M.F. Atiyah and I.G. Macdonald, {\sl Introduction to Commutative Algebra},
Addison-Wesley Pub. Co., 1969.


\bibitem{GL} J. \'Avila Guzm\'an and J. Lazzarin, {\sl A Morita context related to
finite groups acting partially on a ring}, Algebra Discrete Math. 3
(2009), 49-60.

\bibitem{B} D. Bagio, {\sl Partial actions of inductive groupoids on rings},
to appear in Intern. J. Game Theory Algebra, (2011).

\bibitem{BFP} D. Bagio, D. Fl\^ores and A. Paques, {\sl
Partial actions of ordered groupoids on rings}, J. Algebra Appl. 9
(2010), 501-517.

\bibitem{BLP} D. Bagio, J. Lazzarin and A. Paques, {\sl
Crossed products by twisted partial actions: separability,
semisimplicity and Frobenius properties}, Comm. Algebra 38 (2010),
496-508.

\bibitem{CG} S. Caenepeel and E. De Groot, {\sl Galois theory for
weak Hopf algebras}, Rev. Roumaine Math. Pures Appl. 52 (2007),
151-176.

\bibitem{CHR} S. Chase, D. K. Harrison and A. Rosenberg, {\sl
Galois theory and Galois cohomology of commutative rings}, Mem. AMS
52 (1968), 1-19.

\bibitem{D} M. Dokuchaev, {\sl Partial actions: a survey},
Contemporary Math. 537 (2011), 173-184.

\bibitem{DE} M. Dokuchaev and R. Exel, {\sl Associativity of
crossed products by partial actions, enveloping actions and partial
representations}, Trans. AMS 357 (2005), 1931-1952.


\bibitem{DFP} M. Dokuchaev, M. Ferrero and A. Paques, {\sl Partial
actions and Galois theory}, J. Pure Appl. Algebra 208 (2007), 77-87.

\bibitem{E94} R. Exel, {\sl Circle actions on $C^\star$-algebras,
partial automorphisms and generalized Pimsner-Voiculescu exact
sequences}, J. Funct. Anal. 122 (1994), 361-401.

\bibitem{E98} R. Exel, {\sl Partial actions of groups and actions of semigroups},
Proc. AMS 126 (1998), 3481-3494.

\bibitem{Fl} D. Fl\^ores, {A\c c\~ao de Grup\'oides sobre
\'Algebras: Teoremas de Estrutura}, PhD Thesis, UFRGS, Brazil, 2011.

\bibitem{Gil} N.D. Gilbert, {\sl Actions and expansions of ordered groupoids},
J. Pure Appl. Algebra 198 (2005), 175-195.


\bibitem{L} M.V. Lawson, {\sl Inverse Semigroups. The Theory of
Partial Symmetries}, World Scientific Pub. Co, London, 1998.

\bibitem{McC} K. McClanahan, {\sl K-theory for partial crossed
products by discrete groups}, J. Funct. Anal. 130 (1995), 77-117.


\bibitem{Rw} L. H. Rowen, {\sl Ring Theory - Student edition},
Academic Press, 1991.

\bibitem{W} R. Wisbauer, {\sl Foundations of Module and Ring Theory},
Gordon and Breach Sc. Pub., 1991.





\end{thebibliography}
\end{document}